\documentclass{article}

\usepackage{amsmath}
\usepackage{amsthm}

\usepackage[pdftex]{hyperref}
\usepackage{alltt}

\numberwithin{equation}{section}

\begin{document}

\title{Some Identities involving Two Sets of Basis Vectors and the Metric and Mixed Matrices}

\author{M.J. Kronenburg}
\date{}

\maketitle

\begin{abstract}
Given two sets of basis vectors in $n$-dimensional space,
there exists a relation between their lengths and mutual angles,
expressed as relations between the two metric matrices and the mixed matrix.
In this paper these relations are given, and their consequences for
2-dimensional and 3-dimensional space are investigated,
using a computer algebra program for simplifying expressions.
\end{abstract}

\noindent
\textbf{Keywords}: basis vectors, metric matrix, mixed matrix.\\
\textbf{MSC 2010}: 15A24

\section{Introduction}

The two metric matrices and the mixed matrix describe the lengths of and mutual angles
between two sets of basis vectors in $n$-dimensional space \cite{S93}.
There are relations between the two metric matrices and the mixed matrix
which are valid for any two sets of basis vectors \cite{K04,S93}.
These relations can be used in a computer algebra program
to find the general identities of the mutual vector angles in 2-dimensional and 3-dimensional space.
From these general identities, two special cases \cite{S93}, one where one of the metric matrices is the identity matrix,
and one where the mixed matrix is the identity matrix, can be derived.

\section{The Metric and Mixed Matrices in n-dimensional Space}

Let $i=1\cdots n$ be an index over the dimensions in n-dimensional space,
and let $\{\vec{e}_i\}$ be an orthonormal basis in n-dimensional space,
which means that:
\begin{equation}
 \vec{e}_i\cdot\vec{e}_j = \delta_{ij}
\end{equation}
where $\delta_{ij}$ is the Kronecker delta:
\begin{equation}
 \delta_{ij} =
 \begin{cases}
   1 & \text{if $i=j$} \\
   0 & \text{if $i\neq j$} \\
 \end{cases}
\end{equation}
Let the matrices $A=\{A_{ij}\}$ and $A^*=\{A^*_{ij}\}$ have as columns
the coordinates on this orthonormal basis of two sets of basis vectors
$\{\vec{a}_i\}$ and $\{\vec{a}^*_i\}$, so that \cite{T96}:
\begin{equation}
 A_{ij} = \vec{e}_i\cdot\vec{a}_j
\end{equation}
\begin{equation}
 A^*_{ij} = \vec{e}_i\cdot\vec{a}^*_j
\end{equation}
Let the vector $\vec{r}$ have coordinates on the orthonormal basis $\{\vec{e}_i\}$
and let the vectors $\vec{x}$ and $\vec{x}^*$ be the same as vector $\vec{r}$,
but with coordinates on the basis $\{\vec{a}_i\}$ and $\{\vec{a}^*_i\}$, respectively.
Then:
\begin{equation}
 \vec{r} = A\cdot\vec{x} = A^*\cdot\vec{x}^*
\end{equation}
From this equation follows:
\begin{equation}\label{vecx}
 \vec{x}^* = A^{*-1}\cdot A\cdot\vec{x}
\end{equation}
Let the lengths of and the mutual angles between the $\{\vec{a}_i\}$ and the $\{\vec{a}_i^*\}$
be given, that is the following matrices are given:
\begin{equation}\label{gdef}
 G = A^T\cdot A
\end{equation}
\begin{equation}\label{gstardef}
 G^* = A^{*T}\cdot A^*
\end{equation}
\begin{equation}\label{qdef}
 Q = A^T\cdot A^*
\end{equation}
The matrices $G$ and $G^*$ are the metric matrices \cite{G92,S93} and the matrix $Q$ is the mixed matrix \cite{S93}.
Their elements are $G_{ij}=\vec{a}_i\cdot\vec{a}_j$, $G^*_{ij}=\vec{a}^*_i\cdot\vec{a}^*_j$
and $Q_{ij}=\vec{a}_i\cdot\vec{a}^*_j$.
$G$ and $G^*$ are symmetric, and when $G$ (or $G^*$) is given, a matrix $A$ (or $A^*$) can be computed
with Cholesky decomposition \cite{GL83}.
Furthermore:
\begin{equation}\label{rr}
 \vec{r}\cdot\vec{r} = \vec{x}\cdot G\cdot\vec{x} = \vec{x}^*\cdot G^*\cdot\vec{x}^*
   = \vec{x}\cdot Q\cdot\vec{x}^* = \vec{x}^*\cdot Q^T\cdot\vec{x}
\end{equation}
From (\ref{qdef}) follows:
\begin{equation}\label{qres1}
 A^* = (A^{-1})^T\cdot Q
\end{equation}
and equivalently:
\begin{equation}\label{qres2}
 A = (A^{*-1})^T\cdot Q^T
\end{equation}
Now (\ref{vecx}) can also be written as \cite{S93}:
\begin{equation}
 \vec{x}^* = Q^{-1}\cdot G\cdot\vec{x}
\end{equation}
and as
\begin{equation}
 \vec{x}^* = G^{*-1}\cdot Q^T\cdot\vec{x}
\end{equation}
Combination of these two identities yields \cite{K04}:
\begin{equation}\label{gstarid}
 G^* = Q^T\cdot G^{-1}\cdot Q
\end{equation}
and equivalently
\begin{equation}\label{gstaridinv}
 G = Q\cdot G^{*-1}\cdot Q^T
\end{equation}
These last two identities can be used to derive some general identities
involving the mutual angles between the basis vectors $\{\vec{a}_i\}$ and $\{\vec{a}_i^*\}$
in 2-dimensional and 3-dimensional space.\\
In the following, when $\cos(\phi)=x$ is known, because $\phi$ is an angle between two vectors,
we can take $0\leq\phi\leq\pi$, and $\phi=\arccos(x)$.

\section{The 2-dimensional case}

In the 2-dimensional case, let $\alpha_{12}$ be the angle between $\vec{a}_1$ and $\vec{a}_2$,
and $\beta_{12}$ be the angle between $\vec{a}_1^*$ and $\vec{a}_2^*$,
and let $\gamma_{ij}$ be the angle between $\vec{a}_i$ and $\vec{a}_j^*$.
Then the metric and mixed matrices are:
\begin{equation}
 G =
\left(
\begin{matrix}
 |\vec{a}_1|^2 & |\vec{a}_1||\vec{a}_2|\cos(\alpha_{12}) \\
 |\vec{a}_1||\vec{a}_2|\cos(\alpha_{12}) & |\vec{a}_2|^2
\end{matrix}
\right)
\end{equation}
\begin{equation}\label{g2star}
 G^* =
\left(
\begin{matrix}
 |\vec{a}_1^*|^2 & |\vec{a}_1^*||\vec{a}_2^*|\cos(\beta_{12}) \\
 |\vec{a}_1^*||\vec{a}_2^*|\cos(\beta_{12}) & |\vec{a}_2^*|^2
\end{matrix}
\right)
\end{equation}
\begin{equation}\label{Q2def}
 Q =
\left(
\begin{matrix}
 |\vec{a}_1||\vec{a}_1^*|\cos(\gamma_{11}) & |\vec{a}_1||\vec{a}_2^*|\cos(\gamma_{12}) \\
 |\vec{a}_2||\vec{a}_1^*|\cos(\gamma_{21}) & |\vec{a}_2||\vec{a}_2^*|\cos(\gamma_{22})
\end{matrix}
\right)
\end{equation}
By using a computer algebra program, $G^*=Q^T\cdot G^{-1}\cdot Q$ can be evaluated (see below),
resulting in:
\begin{equation}\label{g2res11}
 \cos^2(\gamma_{11})+\cos^2(\gamma_{21})-2\cos(\alpha_{12})\cos(\gamma_{11})\cos(\gamma_{21})=\sin^2(\alpha_{12})
\end{equation}
\begin{equation}\label{g2res22}
 \cos^2(\gamma_{12})+\cos^2(\gamma_{22})-2\cos(\alpha_{12})\cos(\gamma_{12})\cos(\gamma_{22})=\sin^2(\alpha_{12})
\end{equation}
\begin{equation}\label{g2res12}
 \begin{split}
 \cos(\beta_{12}) = \frac{1}{\sin^2(\alpha_{12})}
   [&\cos(\gamma_{11})\cos(\gamma_{12})+\cos(\gamma_{21})\cos(\gamma_{22}) \\
    &- \cos(\alpha_{12})(\cos(\gamma_{11})\cos(\gamma_{22})+\cos(\gamma_{12})\cos(\gamma_{21}))]
\end{split}
\end{equation}
In this case $\sin(\alpha_{12})$ is the area of a parallelogram with side lengths $1$
and with $\alpha_{12}$ the mutual angle between the sides.\\
Subtracting (\ref{g2res22}) from (\ref{g2res11}) results in:
\begin{equation}
 \cos(\alpha_{12}) = \frac{\cos^2(\gamma_{11})+\cos^2(\gamma_{21})-\cos^2(\gamma_{12})-\cos^2(\gamma_{22})}
    {2(\cos(\gamma_{11})\cos(\gamma_{21})-\cos(\gamma_{12})\cos(\gamma_{22}))}
\end{equation} 
This equation can only be used when the denominator is not zero.
When this denominator is zero, it follows from (\ref{g2res11}) and (\ref{g2res22}) that the numerator is also zero,
so that (\ref{g2res11}) and (\ref{g2res22}) are identical.
Using $\sin^2(\alpha_{12})=1-\cos^2(\alpha_{12})$, this equation is quadratic in $\cos(\alpha_{12})$,
and solving this quadratic equation yields that in this case
$\cos(\alpha_{12})=\cos(\gamma_{11}\pm\gamma_{21})=\cos(\gamma_{22}\pm\gamma_{12})$,
where the correct sign must be chosen.
When this denominator is not zero,
and when the four $\gamma_{ij}$ are given, $\alpha_{12}$ and $\beta_{12}$ can be computed
with the last two identities.

\section{The 3-dimensional case}

In the 3-dimensional case, let $\alpha_{ij}$ be the angle between $\vec{a}_i$ and $\vec{a}_j$,
and $\beta_{ij}$ be the angle between $\vec{a}_i^*$ and $\vec{a}_j^*$,
and let $\gamma_{ij}$ be the angle between $\vec{a}_i$ and $\vec{a}_j^*$.
Then the metric and mixed matrices are:
\begin{equation}
 G =
\left(
\begin{matrix}
 |\vec{a}_1|^2 & |\vec{a}_1||\vec{a}_2|\cos(\alpha_{12}) & |\vec{a}_1||\vec{a}_3|\cos(\alpha_{13}) \\
 |\vec{a}_1||\vec{a}_2|\cos(\alpha_{12}) & |\vec{a}_2|^2 & |\vec{a}_2||\vec{a}_3|\cos(\alpha_{23}) \\
 |\vec{a}_1||\vec{a}_3|\cos(\alpha_{13}) & |\vec{a}_2||\vec{a}_3|\cos(\alpha_{23}) & |\vec{a}_3|^2
\end{matrix}
\right)
\end{equation}
\begin{equation}\label{g3star}
 G^* =
\left(
\begin{matrix}
 |\vec{a}_1^*|^2 & |\vec{a}_1^*||\vec{a}_2^*|\cos(\beta_{12}) & |\vec{a}_1^*||\vec{a}_3^*|\cos(\beta_{13}) \\
 |\vec{a}_1^*||\vec{a}_2^*|\cos(\beta_{12}) & |\vec{a}_2^*|^2 & |\vec{a}_2^*||\vec{a}_3^*|\cos(\beta_{23}) \\
 |\vec{a}_1^*||\vec{a}_3^*|\cos(\beta_{13}) & |\vec{a}_2^*||\vec{a}_3^*|\cos(\beta_{23}) & |\vec{a}_3^*|^2
\end{matrix}
\right)
\end{equation}
\begin{equation}
 Q =
\left(
\begin{matrix}
 |\vec{a}_1||\vec{a}_1^*|\cos(\gamma_{11}) & |\vec{a}_1||\vec{a}_2^*|\cos(\gamma_{12}) & |\vec{a}_1||\vec{a}_3^*|\cos(\gamma_{13}) \\
 |\vec{a}_2||\vec{a}_1^*|\cos(\gamma_{21}) & |\vec{a}_2||\vec{a}_2^*|\cos(\gamma_{22}) & |\vec{a}_2||\vec{a}_3^*|\cos(\gamma_{23}) \\
 |\vec{a}_3||\vec{a}_1^*|\cos(\gamma_{31}) & |\vec{a}_3||\vec{a}_2^*|\cos(\gamma_{32}) & |\vec{a}_3||\vec{a}_3^*|\cos(\gamma_{33})
\end{matrix}
\right)
\end{equation}
By using a computer algebra program, $G^*=Q^T\cdot G^{-1}\cdot Q$ can be evaluated (see below),
resulting in:
\begin{equation}
 \Delta = 2\cos(\alpha_{12})\cos(\alpha_{13})\cos(\alpha_{23})-\cos^2(\alpha_{12})-\cos^2(\alpha_{13})-\cos^2(\alpha_{23})+1
\end{equation}
\begin{equation}
 \Omega_1 = \cos(\alpha_{12})\cos(\alpha_{13})-\cos(\alpha_{23})
\end{equation}
\begin{equation}
 \Omega_2 = \cos(\alpha_{12})\cos(\alpha_{23})-\cos(\alpha_{13})
\end{equation}
\begin{equation}
 \Omega_3 = \cos(\alpha_{13})\cos(\alpha_{23})-\cos(\alpha_{12})
\end{equation}
\begin{equation}\label{g3res11}
\begin{split}
   2[ & \cos(\gamma_{21})\cos(\gamma_{31})\Omega_1
   + \cos(\gamma_{11})\cos(\gamma_{31})\Omega_2 + \cos(\gamma_{11})\cos(\gamma_{21})\Omega_3 ] \\
   +&\cos^2(\gamma_{31})\sin^2(\alpha_{12})+\cos^2(\gamma_{21})\sin^2(\alpha_{13})+\cos^2(\gamma_{11})\sin^2(\alpha_{23}) = \Delta
\end{split}
\end{equation}
\begin{equation}\label{g3res22}
\begin{split}
   2[ & \cos(\gamma_{22})\cos(\gamma_{32})\Omega_1
   + \cos(\gamma_{12})\cos(\gamma_{32})\Omega_2 + \cos(\gamma_{12})\cos(\gamma_{22})\Omega_3 ] \\
   +&\cos^2(\gamma_{32})\sin^2(\alpha_{12})+\cos^2(\gamma_{22})\sin^2(\alpha_{13})+\cos^2(\gamma_{12})\sin^2(\alpha_{23}) = \Delta
\end{split}
\end{equation}
\begin{equation}\label{g3res33}
\begin{split}
   2[ & \cos(\gamma_{23})\cos(\gamma_{33})\Omega_1
   + \cos(\gamma_{13})\cos(\gamma_{33})\Omega_2 + \cos(\gamma_{13})\cos(\gamma_{23})\Omega_3 ] \\
   +&\cos^2(\gamma_{33})\sin^2(\alpha_{12})+\cos^2(\gamma_{23})\sin^2(\alpha_{13})+\cos^2(\gamma_{13})\sin^2(\alpha_{23}) = \Delta
\end{split}
\end{equation}
\begin{equation}\label{g3res12}
\begin{split}
 \cos(\beta_{12}) = \frac{1}{\Delta}[&(\cos(\gamma_{21})\cos(\gamma_{32})+\cos(\gamma_{22})\cos(\gamma_{31}))\Omega_1 \\
  +&(\cos(\gamma_{11})\cos(\gamma_{32})+\cos(\gamma_{12})\cos(\gamma_{31}))\Omega_2 \\
  +&(\cos(\gamma_{11})\cos(\gamma_{22})+\cos(\gamma_{12})\cos(\gamma_{21}))\Omega_3 \\
  +&\cos(\gamma_{31})\cos(\gamma_{32})\sin^2(\alpha_{12})+\cos(\gamma_{21})\cos(\gamma_{22})\sin^2(\alpha_{13}) \\
  +&\cos(\gamma_{11})\cos(\gamma_{12})\sin^2(\alpha_{23})]
\end{split}
\end{equation}
\begin{equation}\label{g3res13}
\begin{split}
 \cos(\beta_{13}) = \frac{1}{\Delta}[&(\cos(\gamma_{21})\cos(\gamma_{33})+\cos(\gamma_{23})\cos(\gamma_{31}))\Omega_1 \\
  +&(\cos(\gamma_{11})\cos(\gamma_{33})+\cos(\gamma_{13})\cos(\gamma_{31}))\Omega_2 \\
  +&(\cos(\gamma_{11})\cos(\gamma_{23})+\cos(\gamma_{13})\cos(\gamma_{21}))\Omega_3 \\
  +&\cos(\gamma_{31})\cos(\gamma_{33})\sin^2(\alpha_{12})+\cos(\gamma_{21})\cos(\gamma_{23})\sin^2(\alpha_{13}) \\
  +&\cos(\gamma_{11})\cos(\gamma_{13})\sin^2(\alpha_{23})]
\end{split}
\end{equation}
\begin{equation}\label{g3res23}
\begin{split}
 \cos(\beta_{23}) = \frac{1}{\Delta}[&(\cos(\gamma_{22})\cos(\gamma_{33})+\cos(\gamma_{23})\cos(\gamma_{32}))\Omega_1 \\
  +&(\cos(\gamma_{12})\cos(\gamma_{33})+\cos(\gamma_{13})\cos(\gamma_{32}))\Omega_2 \\
  +&(\cos(\gamma_{12})\cos(\gamma_{23})+\cos(\gamma_{13})\cos(\gamma_{22}))\Omega_3 \\
  +&\cos(\gamma_{32})\cos(\gamma_{33})\sin^2(\alpha_{12})+\cos(\gamma_{22})\cos(\gamma_{23})\sin^2(\alpha_{13}) \\
  +&\cos(\gamma_{12})\cos(\gamma_{13})\sin^2(\alpha_{23})]
\end{split}
\end{equation}
In this case $\sqrt{\Delta}$ is the volume of a parallelepiped with edge lengths $1$
and with the $\alpha_{ij}$ the mutual angles between the edges.

\section{The Special Case G=I}

When $G=I$, or in other words $\vec{a}_i\cdot\vec{a}_j=\delta_{ij}$,
then the basis vectors $\{\vec{a}_i\}$ are an orthonormal basis \cite{S93},
and (\ref{gstarid}) becomes $G^*=Q^T\cdot Q$.
Because $\alpha_{ij}=\pi/2$,\linebreak $\cos(\alpha_{ij})=0$ and $\sin(\alpha_{ij})=1$.

\subsection{The 2-dimensional case}

In the 2-dimensional case, from (\ref{g2res11}), (\ref{g2res22}) and (\ref{g2res12}) follows:
\begin{equation}\label{g1id1}
  \cos^2(\gamma_{11})+\cos^2(\gamma_{21}) = 1
\end{equation}
\begin{equation}\label{g1id2}
  \cos^2(\gamma_{12})+\cos^2(\gamma_{22}) = 1
\end{equation}
\begin{equation}\label{g1id3}
 \cos(\beta_{12}) = \cos(\gamma_{11})\cos(\gamma_{12}) + \cos(\gamma_{21})\cos(\gamma_{22})
\end{equation}
Identities (\ref{g1id1}) and (\ref{g1id2}) also follow from trigonometry.
Using $\cos^2(x)=1-\sin^2(x)$ and $\cos^2(x)-\sin^2(y) = \cos(x+y)\cos(x-y)$ \cite{AS72},
these two identities become:
\begin{equation}\label{g2i1}
 \cos(\gamma_{11}+\gamma_{21})\cos(\gamma_{11}-\gamma_{21}) = 0
\end{equation}
\begin{equation}\label{g2i2}
 \cos(\gamma_{12}+\gamma_{22})\cos(\gamma_{12}-\gamma_{22}) = 0
\end{equation}
Because $\alpha_{12}=\pi /2$, there is some sign combination for which:
\begin{equation}
 \pm\gamma_{11}\pm\gamma_{21} = \frac{\pi}{2} \bmod \pi
\end{equation}
\begin{equation}
 \pm\gamma_{12}\pm\gamma_{22} = \frac{\pi}{2} \bmod \pi
\end{equation}
and because $\cos(\pi /2\bmod\pi)=0$ the two identities follow.\\
Identity (\ref{g1id3}) also follows from trigonometry.\\
Using $\cos(x)\cos(y) = (\cos(x+y)+\cos(x-y))/2$ \cite{AS72}
this identity becomes:
\begin{equation}
 \cos(\beta_{12}) = \frac{1}{2}[\cos(\gamma_{11}+\gamma_{12})+\cos(\gamma_{11}-\gamma_{12})
  +\cos(\gamma_{22}+\gamma_{21})+\cos(\gamma_{22}-\gamma_{21})]
\end{equation}
There is some sign combination for which:
\begin{equation}
 \pm\gamma_{11}\pm\gamma_{12} = \beta_{12}\bmod 2\pi
\end{equation}
\begin{equation}
 \pm\gamma_{21}\pm\gamma_{22} = \beta_{12}\bmod 2\pi
\end{equation}
The other two terms are $\cos(x)+\cos(y)=2\cos((x+y)/2)\cos((x-y)/2)$ \cite{AS72},
and one of these arguments becomes $\pi/2\bmod\pi$.

\subsection{The 3-dimensional case}

In the 3-dimensional case, because $\Delta=1$, $\Omega_i=0$ and $\sin(\alpha_{ij})=1$,
it follows from (\ref{g3res11}), (\ref{g3res22}) and (\ref{g3res33}):
\begin{equation}
 \cos^2(\gamma_{11})+\cos^2(\gamma_{21})+\cos^2(\gamma_{31}) = 1
\end{equation} 
\begin{equation}
 \cos^2(\gamma_{12})+\cos^2(\gamma_{22})+\cos^2(\gamma_{32}) = 1
\end{equation} 
\begin{equation}
 \cos^2(\gamma_{13})+\cos^2(\gamma_{23})+\cos^2(\gamma_{33}) = 1
\end{equation} 
and it follows from (\ref{g3res12}), (\ref{g3res13}) and (\ref{g3res23}):
\begin{equation}
 \cos(\beta_{12}) = \cos(\gamma_{11})\cos(\gamma_{12})+\cos(\gamma_{21})\cos(\gamma_{22})
   +\cos(\gamma_{31})\cos(\gamma_{32})
\end{equation}
\begin{equation}
 \cos(\beta_{13}) = \cos(\gamma_{11})\cos(\gamma_{13})+\cos(\gamma_{21})\cos(\gamma_{23})
   +\cos(\gamma_{31})\cos(\gamma_{33})
\end{equation}
\begin{equation}
 \cos(\beta_{23}) = \cos(\gamma_{12})\cos(\gamma_{13})+\cos(\gamma_{22})\cos(\gamma_{23})
   +\cos(\gamma_{32})\cos(\gamma_{33})
\end{equation}

\section{The Special Case Q=I}

When $Q=I$, or in other words $\vec{a}_i\cdot\vec{a}_j^*=\delta_{ij}$,
then the basis vectors $\{\vec{a}_i\}$ and $\{\vec{a}_i^*\}$ are
mutually reciprocal \cite{S93},
and (\ref{gstarid}) becomes $G^*=G^{-1}$.
Then it follows from (\ref{qres1}) and (\ref{qres2}):
\begin{equation}
 A^* = (A^{-1})^T
\end{equation}
and equivalently
\begin{equation}
 A = (A^{*-1})^T
\end{equation}
Furthermore from (\ref{rr}):
\begin{equation}
  \vec{r}\cdot\vec{r} = \vec{x}\cdot\vec{x}^*
\end{equation}
From the definition of $Q$ it also follows in this case that the only $\cos(\gamma_{ij})$
unequal to zero are when $i=j$, and\begin{equation}\label{Qident}
 |\vec{a}_i||\vec{a}_i^*|\cos(\gamma_{ii}) = 1
\end{equation}
from which follows that:
\begin{equation}\label{gammapos}
  \cos(\gamma_{ii})>0
\end{equation}
and by taking $0<\alpha_{ij}<\pi$ we also have:
\begin{equation}\label{alphapos}
  \sin(\alpha_{ij})>0
\end{equation}

\subsection{The 2-dimensional case}

In the 2-dimensional case, (\ref{g2res11}) and (\ref{g2res22}) with (\ref{gammapos}) and (\ref{alphapos}) give:
\begin{equation}
  \cos(\gamma_{11}) = \sin(\alpha_{12})
\end{equation}
\begin{equation}
  \cos(\gamma_{22}) = \sin(\alpha_{12})
\end{equation}
Combination with (\ref{Qident}) yields:
\begin{equation}
 |\vec{a}_1^*| = \frac{1}{|\vec{a}_1|\sin(\alpha_{12})}
\end{equation}
\begin{equation}
 |\vec{a}_2^*| = \frac{1}{|\vec{a}_2|\sin(\alpha_{12})}
\end{equation}
and with (\ref{g2res12}):
\begin{equation}
  \cos(\beta_{12}) = -\cos(\alpha_{12})
\end{equation}

\subsection{The 3-dimensional case}

In the 3-dimensional case, (\ref{g3res11}), (\ref{g3res22}) and (\ref{g3res33}) with (\ref{gammapos}) and (\ref{alphapos}) give:
\begin{equation}
 \cos(\gamma_{11}) \sin(\alpha_{23}) = \sqrt{\Delta}
\end{equation}
\begin{equation}
 \cos(\gamma_{22}) \sin(\alpha_{13}) = \sqrt{\Delta}
\end{equation}
\begin{equation}
 \cos(\gamma_{33}) \sin(\alpha_{12}) = \sqrt{\Delta}
\end{equation}
Combination with (\ref{Qident}) yields:
\begin{equation}
  |\vec{a}_1^*| = \frac{\sin(\alpha_{23})}{|\vec{a}_1|\sqrt{\Delta}}
\end{equation}
\begin{equation}
  |\vec{a}_2^*| = \frac{\sin(\alpha_{13})}{|\vec{a}_2|\sqrt{\Delta}}
\end{equation}
\begin{equation}
  |\vec{a}_3^*| = \frac{\sin(\alpha_{12})}{|\vec{a}_3|\sqrt{\Delta}}
\end{equation}
and with (\ref{g3res12}), (\ref{g3res13}) and (\ref{g3res23}):
\begin{equation}
 \cos(\beta_{12}) = \frac{\cos(\alpha_{13})\cos(\alpha_{23})-\cos(\alpha_{12})}
  {\sin(\alpha_{13})\sin(\alpha_{23})}
\end{equation}
\begin{equation}
 \cos(\beta_{13}) = \frac{\cos(\alpha_{12})\cos(\alpha_{23})-\cos(\alpha_{13})}
  {\sin(\alpha_{12})\sin(\alpha_{23})}
\end{equation}
\begin{equation}
 \cos(\beta_{23}) = \frac{\cos(\alpha_{12})\cos(\alpha_{13})-\cos(\alpha_{23})}
  {\sin(\alpha_{12})\sin(\alpha_{13})}
\end{equation}

\section{The Inverse Identities}

The matrix identity (\ref{gstarid}), that is $G^*=Q^T\cdot G^{-1}\cdot Q$,
was used for deriving the above identities.
When using matrix identity (\ref{gstaridinv}), that is $G=Q\cdot G^{*-1}\cdot Q^T$,
then the resulting identities are obtained from the above identities
by replacing $|\vec{a}_i|$ by $|\vec{a}^*_i|$, $|\vec{a}^*_i|$ by $|\vec{a}_i|$,
$\alpha_{ij}$ by $\beta_{ij}$, $\beta_{ij}$ by $\alpha_{ij}$,
and $\gamma_{ij}$ by $\gamma_{ji}$.

\section{Computer Program}

The Mathematica$^{\textregistered}$ \cite{W03} programs used to compute the expressions
are given below.\\
The 2-dimensional case:
\begin{alltt}
G:=\{\{r1^2,r1 r2 Cos[alpha12]\},
    \{r1 r2 Cos[alpha12],r2^2\}\}
Q:=\{\{r1 s1 Cos[gamma11],r1 s2 Cos[gamma12]\},
    \{r2 s1 Cos[gamma21],r2 s2 Cos[gamma22]\}\}
FullSimplify[Transpose[Q].Inverse[G].Q]
\end{alltt}
The 3-dimensional case:
\begin{alltt}
G:=\{\{r1^2,r1 r2 Cos[alpha12],r1 r3 Cos[alpha13]\},
    \{r1 r2 Cos[alpha12],r2^2,r2 r3 Cos[alpha23]\},
    \{r1 r3 Cos[alpha13],r2 r3 Cos[alpha23],r3^2\}\}
Q:=\{\{r1 s1 Cos[gamma11],r1 s2 Cos[gamma12],r1 s3 Cos[gamma13]\},
    \{r2 s1 Cos[gamma21],r2 s2 Cos[gamma22],r2 s3 Cos[gamma23]\},
    \{r3 s1 Cos[gamma31],r3 s2 Cos[gamma32],r3 s3 Cos[gamma33]\}\}
FullSimplify[Transpose[Q].Inverse[G].Q]
\end{alltt}

\pdfbookmark[0]{References}{}

\end{document}